\documentclass[twocolumn]{autart}    

\usepackage{graphicx,import}          
\usepackage{color}
\usepackage{amsmath,amssymb,amsfonts}                             
\everymath{\displaystyle}
\usepackage{hyperref}
\usepackage[authoryear]{natbib}
\usepackage{epstopdf}

\newtheorem{assumption}{Assumption}

\newtheorem{lemma}{Lemma}
\newtheorem{remark}{Remark}
\newtheorem{problem}{Problem}

\begin{document}

\begin{frontmatter}

\title{Constrained Stabilization on the $n-$Sphere\thanksref{footnoteinfo}} 

\thanks[footnoteinfo]{This research work is partially supported by NSERC-DG RGPIN-2020-04759, the European Research Council (ERC), the EU H2020 Co4Robots, the SSF COIN project, the Swedish Research Council (VR) and the Knut och Alice Wallenberg Foundation.}

\author[Soulaimane]{Soulaimane Berkane}\ead{soulaimane.berkane@uqo.ca},               
\author[Dimos]{Dimos V. Dimarogonas}\ead{dimos@kth.se}  

\address[Soulaimane]{ Department of Computer Science and Engineering, University   of   Quebec   in   Outaouais, 101 St-Jean Bosco, Gatineau, QC, J8X 3X7, Canada}
\address[Dimos]{Division of Decision and Control Systems, School of Electrical Engineering and Computer Science, KTH Royal Institute of Technology, SE-100 44, Stockholm, Sweden}  

\begin{keyword}                           
Constrained control; unit $n-$sphere; conic constraints, obstacle avoidance, Euclidean sphere world.
\end{keyword}

\begin{abstract}
    We solve the stabilization problem on the $n-$sphere in the presence of conic constraints. We use the stereographic projection to map this problem to the classical navigation problem on $\mathbb{R}^n$ in the presence of spherical obstacles. As a consequence, any obstacle avoidance algorithm for navigation in the Euclidean space can be used to solve the given problem on the $n-$sphere. We illustrate the effectiveness of the approach using the kinematics of the reduced attitude model on the $2-$sphere.
\end{abstract}
\end{frontmatter}

\section{Introduction}
Different mechanical systems of interest have state components that are restricted to evolve on the $n-$sphere.  Example of such systems are the spherical pendulum \citep{Shiriaev2004SwingingIntegrals}, the nonholonomic rolling sphere \citep{Das2004ExponentialSphere}, the reduced-attitude or spin-axis stabilization of rigid bodies \citep{Bullo1995ControlStabilization,Tsiotras1994Spin-axisTorques} and the thrust-vector control for quad-rotor aircraft \citep{Hua2009ADrones}. \cite{Brockett1973LieSpheres} developed a  theory for the most elementary class of control problems defined on spheres where he discussed issues related to controllability, observability and optimal control without an explicit search for control laws. \cite{Bullo1995ControlStabilization} proposed a geometric approach to design controllers for control systems on the sphere relying on the notion of geodesics. The $n-$sphere is not diffeomorphic to a Euclidean space (it is a compact manifold without boundary) and, hence, there exist no smooth control law that globally stabilizes an equilibrium point on the $n-$sphere since the domain of asymptotic stability of any critical point of a continuous vector field needs to be diffeomorphic to a Euclidean space. Recently, hybrid approaches have been proposed to guarantee global asymptotic and exponential stabilization on the $n-$sphere \citep{Mayhew2013GlobalBodies,Casau2015GlobalSphere} and also on the group of rotations $\mathbb{SO}(3)$ \citep{Berkane2017g,Berkane2017f}. 

However, in contrast to the unconstrained stabilization problem, there are only few research works that have considered the constrained stabilization problem on the $n-$sphere. This problem is relevant in different applications such as the pointing maneuver of a space telescope towards a given target ({\it e.g.,} planets and galaxies) during which the telescope's sensitive optical instruments must not be pointed towards bright objects such as stars. The control barrier functions (CBFs) approach on manifolds \citep{Wu2015Safety-criticalManifolds} can be used to solve the constrained stabilization problem on the $n-$sphere. However, this comes at the expense of solving a state-dependent online quadratic program and, besides, the domain of attraction is not characterized. For the particular case of the $2-$sphere, the constrained stabilization problem can be lifted to the constrained (full) attitude stabilization problem where different approaches exist. \cite{Spindler2002AttitudeDirection} proposed a geometric control law that minimizes a given cost function to solve the problem of maneuvering a rigid spacecraft attitude from rest to rest while avoiding a single forbidden direction. In \citep{Lee2014FeedbackPotentials}, a logarithmic barrier potential function is used to synthesize a quaternion-based feedback controller that solves the attitude reorientation of a rigid body spacecraft in the presence of multiple attitude-constrained zones. Another potential-based approach for the constrained attitude control on $\mathbb{SO}(3)$ has been proposed in \citep{Kulumani2017ConstrainedSO3}. 

In this work we consider the constrained stabilization problem of dynamical systems evolving on a configuration space defined by the unit $n-$sphere. The considered forbidden zones are conic-type constraints in the sense that we force the state on the $n-$sphere (which can be seen as a unit axis) to keep minimum safety angles with respect to some given set of other unit axes. Our proposed solution consists first in showing that the considered constrained $n-$sphere manifold is diffeomorphic, via the stereographic projection \citep{Helmke1996OptimizationSystems}, to a Euclidean space punctured by spherical obstacles. Then, by considering a generic driftless system on the $n-$sphere, we prove that the pushforward vector field in the new stereographic coordinates is feedback linearizable. Therefore, we are able to map the given constrained stabilization problem on the $n-$sphere to the well-established and treated obstacle avoidance problem in $\mathbb{R}^n$ which allows us to benefit from the many studies in the latter field. For instance, one can use navigation functions \citep{Koditschek1990RobotBoundary} to obtain an almost global result or even global results with hybrid control techniques \citep{Berkane2019}. We show that the qualitative properties (e.g., invariance, stability, region of attraction) of any static obstacle avoidance controller on $\mathbb{R}^n$ are preserved for the resulting safety controller on the $n-$sphere. {\bf Notation}: We use $\mathbb{N}$, $\mathbb{R}$ and $\mathbb{R}_{\geq 0}$ to denote, respectively, the sets of positive integers, real and non-negative real numbers. $\mathbb{R}^n$ denotes the $n-$dimensional Euclidean space. $I_n$ denotes the $n\times n$ identity matrix and $e_k$ corresponds to the $k-$th column of $I_n$. The Euclidean norm of $x\in\mathbb{R}^n$ is defined as $\|x\|=\sqrt{x^\top x}$ where $(\cdot)^{\top}$ denotes the transpose of $(\cdot)$. The topological interior (resp. boundary) of a subset $\mathcal{S}$ of a metric space is denoted by $\mathbf{int}(S)$ (resp. $\partial\mathcal{S}$). For a multi-variable function $f(x_1,\cdots,x_n)$, we denote by $\nabla_if$ the gradient of $f$ with respect to the $i-$th argument $x_i$.
\vspace{-0.45cm}
\section{Problem Formulation}\label{section:problem}
\vspace{-0.3cm}
The unit $n-$sphere is an $n-$dimensional manifold that is embedded in the Euclidean space $\mathbb{R}^{n+1}$ and defined as
$
    \mathbb{S}^n:=\{x\in\mathbb{R}^{n+1}: \|x\|=1\}.
$
The tangent space to $\mathbb{S}^n$ at a given point $x$ is defined by the $n-$dimensional hyperplane
$
    \mathsf{T}_x\mathbb{S}^n=\{z\in\mathbb{R}^{n+1}: z^\top x=0\},
$
which represents all vectors in $\mathbb{R}^{n+1}$ that are perpendicular to $x\in\mathbb{S}^n$. $\mathbb{S}^n$ is a metric space if we pair it with the geodesic distance $\mathbf{d}(x,y):=\arccos(x^\top y),\forall x,y\in\mathbb{S}^n.
$
We consider the following driftless system on $\mathbb{S}^n$:
\begin{align}\label{eq:dx}
    \dot x=\Pi(x)u,
\end{align}
where $u\in\mathbb{R}^m$ is the control input and $\Pi:\mathbb{R}^{n+1}\to\mathbb{R}^{(n+1)\times m}$ is a smooth matrix-valued function such that $\mathrm{Im}(\Pi(x))\subseteq\mathsf{T}_x\mathbb{S}^n$. The condition $\mathrm{Im}(\Pi(x))\subseteq\mathsf{T}_x\mathbb{S}^n$ implies that $\dot x\in\mathsf{T}_x\mathbb{S}^n$ or $x^\top\dot x=0$ which guarantees forward invariance of $\mathbb{S}^n$ under the dynamics \eqref{eq:dx} since $\|x\|^2$ remains constant regardless of the input $u$. Our goal is to propose a constrained stabilization strategy on $\mathbb{S}^n$ in the presence of the following $(I+1)$ conic constraints:
\begin{align}\label{eq:Oi}
    \mathcal{O}_i=\{x\in\mathbb{S}^n: x^\top a_i>\cos(\theta_i)\}, i\in\mathbb{I}=\{0,\cdots,I\}
\end{align}
where $a_i\in\mathbb{S}^n$ is the center of $\mathcal{O}_i$ and $\theta_i\in(0,\pi/2)$ is the smallest angle (between $x$ and $a_i$) allowed in the free region. We define our {\it constrained space on $\mathbb{S}^n$} as $\mathcal{M}:=\mathbb{S}^n\setminus\cup_{i\in\mathbb{I}}\mathcal{O}_i=
  \{x\in\mathbb{S}^n:\mathbf{d}(x,a_i)\geq\theta_i,\;\forall i\in\mathbb{I}\}$.
\begin{assumption}\label{assumption:general}
The following assumptions hold:
\begin{enumerate}
    \item For all $x\in\mathbb{S}^n\setminus\{e_{n+1}\}$, $\mathrm{rank}(\Pi(x))=n$.
    \item For all $i,i^\prime\in\mathbb{I}$ with $i\neq i^\prime$, $a_i^\top a_{i^\prime}<\cos(\theta_i+\theta_{i^\prime})$.
    \item $a_0=e_{n+1}$.
    \item $x(0)\in\mathcal{M}$ and $x_d\in\mathbf{int}(\mathcal{M})$.
\end{enumerate}
\end{assumption}
Item $(1)$ of Assumption \ref{assumption:general} is a controlability assumption that imposes the fact that we can steer any point on $\mathbb{S}^n\setminus\{e_{n+1}\}$ in any direction by appropriately choosing the control input $u$. 
Item $(2)$ imposes that the closures of the constraint zones $\mathcal{O}_i$ are pairwise disjoint. In item $(3)$ we assume, without loss of generality, that the obstacle $\mathcal{O}_0$'s axis coincides with the coordinate axis $e_{n+1}$. Finally, item $(4)$ imposes that the initial condition $x(0)$ and the desired reference point $x_d$ must lie in the free space $\mathcal{M}$ and the interior of $\mathcal{M}$, respectively.
\begin{problem}\label{problem1}
For system \eqref{eq:dx} and under Assumption \ref{assumption:general}, design a control law $u=\kappa(x,x_d)$ that renders the constrained space $\mathcal{M}$ forward invariant and the target point $x=x_d$ an asymptotically (or exponentially) stable equilibrium with a region of attraction $\mathcal{R}(x_d)\subseteq\mathcal{M}$.   
\end{problem}
%
\vspace{-0.35cm}
\section{Main Results}\label{section:main}
\begin{figure}
    \centering
    \includegraphics[width=0.8\columnwidth]{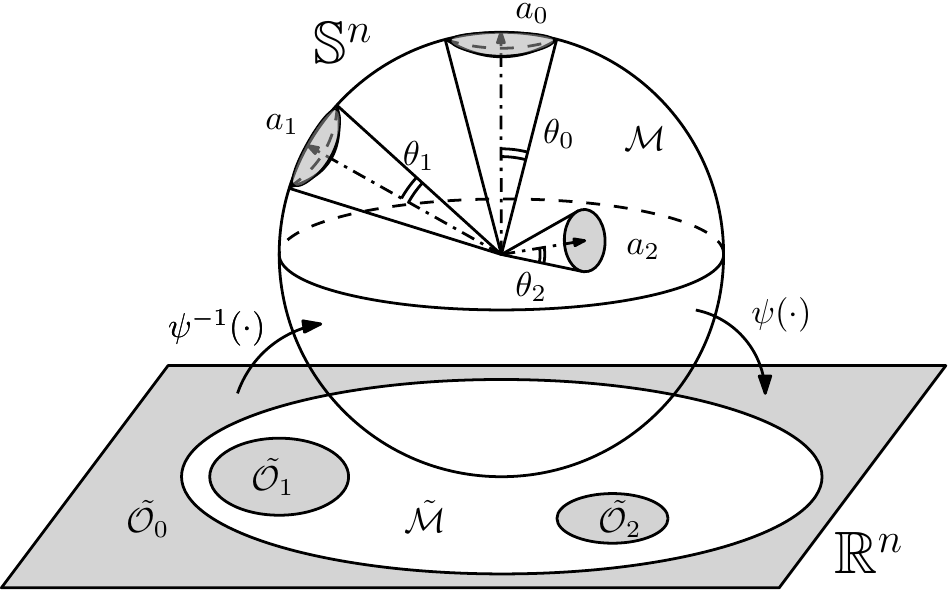}
    \caption{The stereographic projection $\psi(\cdot)$  maps the constrained space $\mathcal{M}$ to the Euclidean sphere world $\tilde{\mathcal{M}}$.}
    \label{fig:sphere_map}
\end{figure}

\vspace{-0.25cm}
The {\it stereographic projection} is defined by the map $\psi:\mathbb{S}^n\setminus\{e_{n+1}\}\to\mathbb{R}^n$ \citep{Helmke1996OptimizationSystems}
\begin{align}
    \label{eq:psi}
    \psi(x)=\frac{\begin{bmatrix}
I_{n}&0_{n\times 1}
\end{bmatrix}x}{1-e_{n+1}^\top x}=:\frac{J_nx}{1-e_{n+1}^\top x}.
\end{align}
Geometrically speaking, the stereographic projection of a point $x\in\mathbb{S}^n\setminus\{e_{n+1}\}$ represents the unique point $\psi(x)$  describing the intersection of the line, that passes by $e_{n+1}$ and $x$, with the hyperplane $\{x\in\mathbb{R}^{n+1}:x_{n+1}=0\}$. The following are some useful properties of this map. 
\begin{lemma}\label{lemma:psi}
The stereographic projection map $\psi$ satisfies:
\begin{itemize}
    \item $\psi$ is a diffeomorphism with the inverse given explicitly by the map $\psi^{-1}:\mathbb{R}^n\to\mathbb{S}^n\setminus\{e_{n+1}\}$ such that
\begin{align}\label{eq:inverse}
    \psi^{-1}(\xi):=\frac{2J_n^\top\xi+(\|\xi\|^2-1)e_{n+1}}{1+\|\xi\|^2}.
    \end{align}
    \item The Jacobian matrix of $\psi(x)$ is given by
\begin{equation}\label{eq:nabla-psi}
    \nabla\psi(x)=\frac{J_n\big((1-e_{n+1}^\top x)I_{n+1}+xe_{n+1}^\top\big)}{(1-e_{n+1}^\top x)^2}.
\end{equation} 
\item For all $x_1,x_2\in\mathbb{S}^n\setminus\{e_{n+1}\}$
\begin{equation}\label{eq:norm-psi}
    \|\psi(x_1)-\psi(x_2)\|^2=\frac{2(1-x_1^\top x_2)}{(1-e_{n+1}^\top x_1)(1-e_{n+1}^\top x_2)}.
\end{equation}
\item If item $(1)$ of Assumption \ref{assumption:general} holds, then $\Sigma(x):=\nabla\psi(x)\Pi(x)$ is full row rank for all $x\in\mathbb{S}^n\setminus\{e_{n+1}\}$.
\end{itemize}
\end{lemma}
\vspace{-0.5cm}
\begin{pf}
The explicit expression $\psi^{-1}$ is taken from \cite[Appendix C.4]{Helmke1996OptimizationSystems}. Also, both $\psi$ and $\psi^{-1}$ are differentiable on their domains of definition and, hence, $\psi$ is a diffeomorphism. The Jacobian of $\psi$ is obtained by direct differentiation of \eqref{eq:psi}. Making use of $ J_n^\top J_n=I_{n+1}-e_{n+1}e_{n+1}^\top$, we have
\begin{equation}\label{eq:psi-norm}
        \|\psi(x)\|^2=\frac{1-(e_{n+1}^\top x)^2}{(1-e_{n+1}^\top x)^2}=\frac{1+e_{n+1}^\top x}{1-e_{n+1}^\top x}.
\end{equation}%
Therefore, for all $x_1,x_2\in\mathbb{S}^n\setminus\{e_{n+1}\}$\vskip -1cm
{\small\begin{align*}
    &\|\psi(x_1)-\psi(x_2)\|^2=\|\psi(x_1)\|^2+\|\psi(x_2)\|^2-2\psi(x_1)^\top\psi(x_2)\\
    &\overset{\eqref{eq:psi-norm}}{=}\frac{1+e_{n+1}^\top x_1}{1-e_{n+1}^\top x_1}+\frac{1+e_{n+1}^\top x_2}{1-e_{n+1}^\top x_2}-\frac{2x_1^\top J_n^\top J_nx_2}{(1-e_{n+1}^\top x_1)(1-e_{n+1}^\top x_2)}\\
    &=\frac{2(1-x_1^\top x_2)}{(1-e_{n+1}^\top x_1)(1-e_{n+1}^\top x_2)}.
\end{align*}}%
Since $\psi$ is a diffeomorphism, we have $\forall x\in\mathbb{S}^n\setminus\{e_{n+1}\}$, $\mathrm{rank}(\nabla\psi(x))=n$. Hence, by the rank-nullity theorem, $\mathrm{dim}(\mathrm{ker}(\nabla\psi(x)))=1$. However, $\nabla\psi(x)(x-e_{n+1})=0$ and thus $\mathrm{ker}(\nabla\psi(x))=\{\alpha(x-e_{n+1}):\alpha\in\mathbb{R}\}$. On the other hand, since $x\in\mathbb{S}^n\setminus\{e_{n+1}\}$, we have $x^\top(x-e_{n+1})=1-x^\top e_{n+1}\neq 0$ which implies that $\alpha(x-e_{n+1})\notin\mathsf{T}_x\mathbb{S}^n$ or $\mathrm{ker}(\nabla\psi(x))\cap\mathrm{Im}(\Pi(x))=\emptyset$. Finally, by applying \cite[Fact 2.10.14., item ii)]{Bernstein2009MatrixFormulas}, we hace
$\mathrm{rank}(\Sigma(x))=\mathrm{rank}(\Pi(x))-\mathrm{dim}(\mathrm{ker}(\nabla\psi(x))\cap\mathrm{Im}(\Pi(x)))=\mathrm{rank}(\Pi(x))=n$.
\end{pf}
We show in the following lemma that $\psi$ maps the constrained space $\mathcal{M}$ to a Euclidean sphere world on $\mathbb{R}^n$ as defined in \citep{Koditschek1990RobotBoundary}; see Fig.~\ref{fig:sphere_map}.
\begin{lemma}\label{lemma:image-of-Oi}
Let Assumption \ref{assumption:general} hold. Then, the image of the constrained zones $\{\mathcal{O}_i\}_{i\in\mathbb{I}}$ via the map $\psi$ is given by the following pairwise disjoint spherical subsets of $\mathbb{R}^n$
\begin{align*}
    &\tilde{\mathcal{O}}_0:=\psi(\mathcal{O}_0)=\{\xi\in\mathbb{R}^n: \|\xi\|>\cot{(\theta_0/2)}\},\\
    &\tilde{\mathcal{O}}_i:=\psi(\mathcal{O}_i)=\{\xi\in\mathbb{R}^n: \|\xi-c_i\|< r_i\}, i\in\mathbb{I}\setminus\{0\},
\end{align*}
where $c_i\in\mathbb{R}^n$ and $r_i\in\mathbb{R}_{>0}$ are given by
\begin{align}
    \label{eq:ci}
    c_i=\frac{J_na_i}{\cos(\theta_i)-a_0^\top a_i},\quad   r_i=\frac{\sin(\theta_i)}{\cos(\theta_i)-a_0^\top a_i},
\end{align}
such that $a_i\in\mathbb{R}^{n+1}$ and $\theta_i\in(0,\pi/2)$ are defined in \eqref{eq:Oi}.
\end{lemma}
\vspace{-0.5cm}
\begin{pf}
First, note that $c_i$ and $r_i$ are well-defined since, by item 2) of Assumption \ref{assumption:general}, it holds that $a_0^\top a_i\leq\cos(\theta_0+\theta_i)<\cos(\theta_i)$ since $\theta_0,\theta_i\in(0,\pi/2)$. Since the function $z\mapsto(1+z)/(1-z)$ is monotonically increasing on $[-1, 1)$, we have that $x^\top a_0>\cos(\theta_0)$ is equivalent to $\|\psi(x)\|^2>(1+\cos(\theta_0))/(1-\cos(\theta_0))=\cot^2(\theta_0/2)$ which proves $\psi(\mathcal{O}_0)$. Now, let $x\in\mathcal{O}_i$ for $i\in\mathbb{I}\setminus\{0\}$:\vskip -1cm
{\small
\begin{equation*}
\begin{aligned}
    &\|\psi(x)-c_i\|^2
    =\|\psi(x)\|^2+\|c_i\|^2-2c_i^\top\psi(x)\\
    &\overset{\eqref{eq:psi-norm}}{=}\frac{1+a_0^\top x}{1-a_0^\top x}+\frac{a_i^\top J_n^\top J_na_i}{(\cos(\theta_i)-a_0^\top a_i)^2}-\frac{2a_i^\top J_n^\top J_nx}{(\cos(\theta_i)-a_0^\top a_i)(1-a_0^\top x)}\\
    &=\frac{1+a_0^\top x}{1-a_0^\top x}+\frac{1-(a_0^\top a_i)^2}{(\cos(\theta_i)-a_0^\top a_i)^2}-\frac{2a_i^\top x-2(a_0^\top a_i)(a_0^\top x)}{(\cos(\theta_i)-a_0^\top a_i)(1-a_0^\top x)}\\
&\overset{\eqref{eq:Oi}}{<}\frac{1+a_0^\top x}{1-a_0^\top x}+\frac{1-(a_0^\top a_i)^2}{(\cos(\theta_i)-a_0^\top a_i)^2}-\frac{2\cos(\theta_i)-2(a_0^\top a_i)(a_0^\top x)}{(\cos(\theta_i)-a_0^\top a_i)(1-a_0^\top x)}\\
    &=\frac{1-\cos(\theta_i)^2}{(\cos(\theta_i)-a_0^\top a_i)^2}=r_i^2
\end{aligned}
\end{equation*}}%
which proves $\psi(\mathcal{O}_i)$. Finally, since $\psi$ is a bijection and $\{\mathcal{O}_i\}_{i\in\mathbb{I}}$ are pairwise disjoint, it follows that the sets $\{\psi(\mathcal{O}_i)\}_{i\in\mathbb{I}}$ are also pairwise disjoint.
\end{pf}
The resulting Euclidean sphere world consists of one large $(n-1)-$sphere $\mathbb{R}^n\setminus\psi(\mathcal{O}_0)$ that bounds the workspace and other smaller disjoint  $(n-1)-$spheres $\psi(\mathcal{O}_i)$ that define obstacles in $\mathbb{R}^n$ that are strictly contained in the interior of the workspace. The obtained {\it constrained space on $\mathbb{R}^n$} is 
$
    \tilde{\mathcal{M}}:=\psi(\mathcal{M})=\mathbb{R}^n\setminus\cup_{i\in\mathbb{I}}\tilde{\mathcal{O}}_i.
$
Now, let us consider the change of variable $\xi:=\psi(x)$. Then, in view of \eqref{eq:dx}, the dynamics of $\xi$ are given by
\begin{equation}\label{eq:dpsi}
\begin{aligned}
    \dot\xi=\nabla\psi(x)\dot x=\nabla\psi(x)\Pi(x)u=:\Sigma(x)u.
\end{aligned}
\end{equation}
Interestingly, according to Lemma \ref{lemma:psi}, $\Sigma(x)$ is full row rank (right invertible) and, hence, its Moore-Penrose pseudo inverse can be explicitly calculated as follows:
\begin{align}\label{eq:Sigma-inverse}
    \Sigma(x)^+=\Sigma(x)^\top(\Sigma(x)\Sigma(x)^\top)^{-1}.
\end{align}
Therefore, by considering a control law of the form
\begin{align}\label{eq:u}
    u=\Sigma(x)^+v,
\end{align}
where $v\in\mathbb{R}^n$ is a virtual control input, one obtains
\begin{align}\label{eq:dxi}
    \dot\xi=v.
\end{align}%
Next, we show in Theorem \ref{theorem} that solving Problem \ref{problem1} on the constrained space $\mathcal{M}$ boils down to solving the following problem on the Euclidean sphere world $\tilde{\mathcal{M}}$.
\begin{problem}\label{problem2}
For system \eqref{eq:dxi} and for a given $\xi_d\in\mathbf{int}(\tilde{\mathcal{M}})$, design a control law $v=\tilde{\kappa}(\xi,\xi_d)$ that renders the constrained space $\tilde{\mathcal{M}}$ forward invariant and the target point $\xi=\xi_d$ an asymptotically (exponentially) stable equilibrium with a region of attraction  $\tilde{\mathcal{R}}(\xi_d)\subseteq\tilde{\mathcal{M}}$.   
\end{problem}

\begin{thm}\label{theorem}
If $v=\tilde{\kappa}(\xi,\xi_d)$ is a control law that solves Problem \ref{problem2} with region of attraction $\tilde{\mathcal{R}}(\xi_d)$ then $u=\Sigma(x)^+\tilde{\kappa}(\psi(x),\psi(x_d))$ is a control law that solves Problem \ref{problem1} with a region of attraction $\mathcal{R}(x_d)=\psi^{-1}(\tilde{\mathcal{R}}(\psi(x_d)))$. Moreover, if $\tilde{\kappa}(\xi,\xi_d)$ is a priori bounded on $\tilde{\mathcal{M}}$ then the control law $u$ is also a priori bounded on $\mathcal{M}$.
\end{thm}
\vspace{-0.3cm}
\begin{pf}
To prove this result we use the comparison theorem \cite[Theorem 3.4.1]{Michel2001QualitativeMappings} to deduce the qualitative properties of all solutions of 
\begin{align}\label{eq:dx-closed}
    \dot x=\Pi(x)\Sigma(x)^+\tilde{\kappa}(\psi(x),\psi(x_d)), x(0)\in\mathcal{M},
\end{align}
from the qualitative properties of solutions of
\begin{align}\label{eq:dxi-closed}
    \dot\xi=\tilde{\kappa}(\xi,\xi_d), \xi(0)=\psi(x(0)), \xi_d=\psi(x_d).
\end{align}
We denote hereafter by $\mathcal{S}_x$ (resp. $\mathcal{S}_\xi$) the set of all solutions to \eqref{eq:dx-closed} (resp. \eqref{eq:dxi-closed}).  It is clear that $\psi(\mathcal{S}_x)\subset\mathcal{S}_\xi$ since, for all $p_x(t,x)\in\mathcal{S}_x$, we have $\dot \psi(p_x)=\nabla\psi(p_x)\dot p_x=\Sigma(p_x)\Sigma(p_x)^+\tilde{\kappa}(\psi(p_x),\psi(x_d))=\tilde{\kappa}(\psi(p_x),\xi_d)$ and hence $\psi(p_x(t,x))\in\mathcal{S}_\xi$. Moreover, in view of \eqref{eq:norm-psi} we have
\begin{align}
     &\|\psi(x)-\psi(x_d)\|^2=\frac{2(1-x^\top x_d)}{(1-e_{n+1}^\top x)(1-e_{n+1}^\top x_d)}\\
     &=\frac{2(1-\cos(\mathbf{d}(x,x_d)))}{(1-\cos(\mathbf{d}(x,e_{n+1})))(1-\cos(\mathbf{d}(x_d,e_{n+1})))}\\
     &=\frac{\sin^2(\mathbf{d}(x,x_d)/2)}{\sin^2(\mathbf{d}(x,e_{n+1})/2)\sin^2(\mathbf{d}(x_d,e_{n+1})/2)}.
\end{align}
Now, using the fact that $x,x_d\in\mathcal{M}\subset\mathbb{R}^n\setminus\mathcal{O}_0$ (i.e., $\mathbf{d}(x,e_{n+1})\geq\theta_0$ and $\mathbf{d}(x_d,e_{n+1})\geq\theta_0$) and the useful identity $(2/\pi)z\leq\sin(z)\leq z$, $z\in[0,\pi/2]$, one deduces 
\begin{align*}
    4\pi^{-2}\mathbf{d}^2(x,x_d)\leq\|\psi(x)-\psi(x_d)\|^2\leq\sin^{-4}(\theta_0)\mathbf{d}^2(x,x_d).
\end{align*}
Since the distance $\|\psi(x)-\psi(x_d)\|$ on the closed set $\tilde{\mathcal{M}}$ is upper and lower bounded by class $\mathcal{K}_\infty$ functions (more precisely positive-power functions) of the distance $\mathbf{d}(x,x_d)$ on the closed set $\mathcal{M}$, we can apply \cite[Theorem 3.4.1]{Michel2001QualitativeMappings} to conclude that forward invariance of $\tilde{\mathcal{M}}$ and $\xi=\xi_d$ (with respect to \eqref{eq:dxi-closed}) implies, respectively, forward invariance of $\mathcal{M}$ and $x=x_d$ (with respect to \eqref{eq:dx-closed}).  Moreover, asymptotic (exponential) stability of the equilibrium $\xi=\xi_d$ implies asymptotic (exponential) stability of the equilibrium $x=x_d$. Now, let $x\in\psi^{-1}(\tilde{\mathcal{R}}(\psi(x_d)))$ and $p_x(t,x)\in\mathcal{S}_x$. Then, one has $\psi(x)\in\tilde{\mathcal{R}}(\psi(x_d))$ and $\psi(p_x(t,x))\in\mathcal{S}_\xi$. It follows that $\lim_{t\to\infty}\psi(p_x(t,x))=\psi(x_d)$ and, by continuity of $\psi^{-1}(\cdot)$, one has $\lim_{t\to\infty}p_x(t,x)=x_d$, which shows that $x\in\mathcal{R}(x_d)$. On the other hand, $x\in\mathcal{R}(x_d)$ implies $\lim_{t\to\infty}p_x(t,x)=x_d$) and, by continuity of $\psi(\cdot)$, $\lim_{t\to\infty}\psi(p_x(t,x))=\psi(x_d)$. This is equivalent to $\psi(x)\in\tilde{\mathcal{R}}(\psi(x_d))$ and thus $x\in\psi^{-1}(\tilde{\mathcal{R}}(\psi(x_d)))$. At the end, since $\Pi(x)$ is smooth, $\nabla\psi(x)$ is continuous on $\mathbb{S}^n\setminus\{e_{n+1}\}$, we have that $\Sigma(x)$ (and hence  $\Sigma(x)^+$) is continuous on $\mathcal{M}$. However, $\mathcal{M}$ is compact (closed and bounded) which implies that the continuous function $\Sigma(x)^+$ must be bounded on $\mathcal{M}$. It follows that $\Sigma(x)^+\tilde{\kappa}(\psi(x),\psi(x_d))$ is a priori bounded control on $\mathcal{M}$ if $\tilde{\kappa}(\xi,\xi_d)$ is a priori bounded on $\tilde{\mathcal{M}}$.
\end{pf}
%
\vspace{-0.5cm}
In the case of a single constraint, one has $\tilde{\mathcal{M}}=\mathbb{R}^n\setminus\tilde{\mathcal{O}}_0$ which, in view of Lemma \ref{lemma:image-of-Oi}, represents the ball bounded by the sphere of radius $\cot(\theta_0/2)$ that is centered at $0$. It is not difficult to show that setting $v=-\gamma(\xi-\xi_d)$ in \eqref{eq:dxi} results in GES of $\xi=\xi_d$ and forward invariance of $\tilde{\mathcal{M}}$. The following corollary follows from Theorem \ref{theorem}. 
\begin{cor}[single constraint]\label{cor:single} 
Consider the kinematics~\eqref{eq:dx} under Assumption \ref{assumption:general} and let $I=0$. Then, the control law $u=-\gamma\Sigma(x)^+(\psi(x)-\psi(x_d)), \gamma>0,$ guarantees global exponential stability of the equilibrium $x=x_d$ and forward invariance of the free space $\mathcal{M}$.
\end{cor}
The global result of Corollary~\ref{cor:single} is related to the well-known Alexandroff one-point compactification in general topology \citep{Alexandroff1924}. In fact, removing a single constraint zone $\mathcal{O}_0$ from the unit $n-$sphere results in a manifold that is diffeomorphic to a Euclidean space and, therefore, global asymptotic stability is possible via a continuous time-invariant feedback.

If we have two or more constraints, the constrained manifold $\tilde{\mathcal{M}}$ is not diffeomorphic to a Euclidean space and, hence, there is a topological obstruction to solve Problem~\ref{problem2} globally with a continuous feedback. Different controllers from the vast literature on obstacle avoidance can be  employed here. For instance,  continuous feedback, {\it e.g.,} \citep{Koditschek1990RobotBoundary, Loizou2017TheTransformation}, can be used to ensure almost global asymptotic stabilization while hybrid feedback, {\it e.g.,} \citep{Berkane2019}, can be used to ensure stronger global asymptotic stabilization. In this work, we consider the navigation functions-based approach of \citep{Koditschek1990RobotBoundary}. In particular, let the navigation function
\begin{align}
    \phi(\xi,\xi_d)=\frac{\|\xi-\xi_d\|^2}{(\|\xi-\xi_d\|^{2k}+\beta(\xi))^{\frac{1}{k}}},\quad k>0,
\end{align}%
where $\beta(\xi)=\Pi_{i\in\mathbb{I}}\beta_i(\xi)$ and
\begin{align}
    &\beta_0(\xi)=\cot^2(\theta_0/2)-\|\xi\|^2,\\
    &\beta_i(\xi)=\|\xi-c_i\|^2-r_i^2,\quad i\in\mathbb{I}\setminus\{0\}.
\end{align}
Note that the parameter $k$ needs to be tuned above a certain threshold in order to eliminate local minima and for $\phi$ to be a valid navigation function. We then consider navigation along the negative gradient of $\phi$ and define 
\begin{align}\label{eq:kappa}
\tilde{\kappa}(\xi,\xi_d):=-\gamma\nabla_1\phi(\xi,\xi_d),\quad\gamma>0.
\end{align}
We state the following corollary that follows  from the result of Theorem \ref{theorem} and \citep{Koditschek1990RobotBoundary}.
\begin{cor}[Two or more constraints]
Consider the kinematics~\eqref{eq:dx} under Assumption \ref{assumption:general}. Then, there exists $K$ such that if $k>K$ the control law $u=\Sigma(x)^+\tilde{\kappa}(\psi(x),\psi(x_d))$, with $\tilde{\kappa}$ defined in \eqref{eq:kappa}, guarantees almost global asymptotic stability of the equilibrium $x=x_d$ and forward invariance of $\mathcal{M}$.
\end{cor}

 \begin{remark}
The feedback linearization approach in \eqref{eq:dpsi}-\eqref{eq:dxi} can be extended to high-order dynamics of the form:
\begin{equation}\label{eq:dx12}
    \begin{aligned}
        \dot x_1&=\Pi(x_1)x_2,\\
        \dot x_q&=x_{q+1}, 2\leq q\leq(l-1),\quad\dot x_l&=u,
    \end{aligned}
\end{equation}
where $x=(x_1,x_2,\cdots,x_l)\in\mathbb{S}^n\times\mathbb{R}^{(l-1)m}$ and $u\in\mathbb{R}^m$. Let the change of variables $\xi=(T_1(x),\cdots,T_l(x))$ such that $T_1(x):=\psi(x_1), T_2(x):=\Sigma(x_1)x_2$ and 
\begin{equation*}
  T_{q+1}(x):=\nabla_1T_q(x)\Pi(x_1)x_2+\textstyle\sum_{p=2}^{q}\nabla_pT_q(x)x_{p+1}
\end{equation*}
for $q=2,\cdots,(l-1)$. By construction $T_q(x)$ depends only on $x_1,\cdots,x_q$ and, hence, $\nabla_{q+1} T_{q+1}(x)=\nabla_q T_q(x)$. It follows that $\nabla_l T_l(x)=\nabla_2 T_2(x)=\Sigma(x_1)$ and the dynamics of the new variables are 
\begin{equation*}
    \begin{aligned}
        \dot \xi_q&=\xi_{q+1},\quad1\leq q\leq l-1,\\
        \dot \xi_l&=\nabla_1T_l(x)\Pi(x_1)x_2+\textstyle\sum_{p=2}^{l-1}\nabla_pT_l(x)x_{p+1}+\Sigma(x_1)u.
    \end{aligned}
\end{equation*}%
Therefore, the control law
\begin{align*}
u=\Sigma(x_1)^+\big(v-\nabla_1T_l(x)\Pi(x_1)x_2-\textstyle\sum_{p=2}^{l-1}\nabla_pT_l(x)x_{p+1}\big),
\end{align*}
where $v\in\mathbb{R}^n$, results in the linear dynamics
\begin{equation}\label{eq:dxi12}
    \begin{aligned}
        \dot \xi_q&=\xi_{q+1},1\leq q\leq l-1,\quad\dot \xi_l&=v.
    \end{aligned}
\end{equation}
\end{remark}

\section{Example}\label{section:example}
\begin{figure}
\includegraphics[width=0.495\columnwidth]{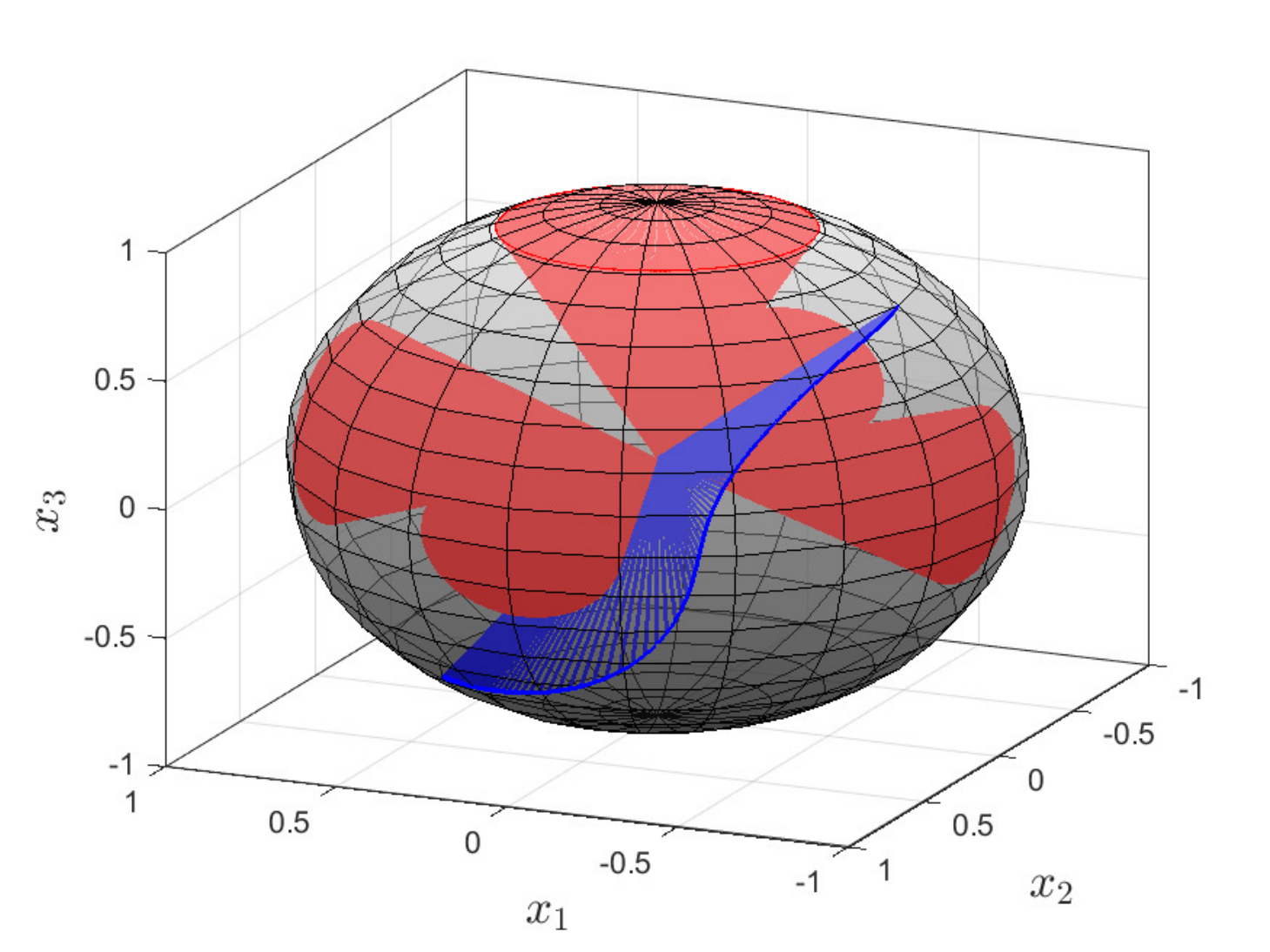}
\includegraphics[width=0.495\columnwidth]{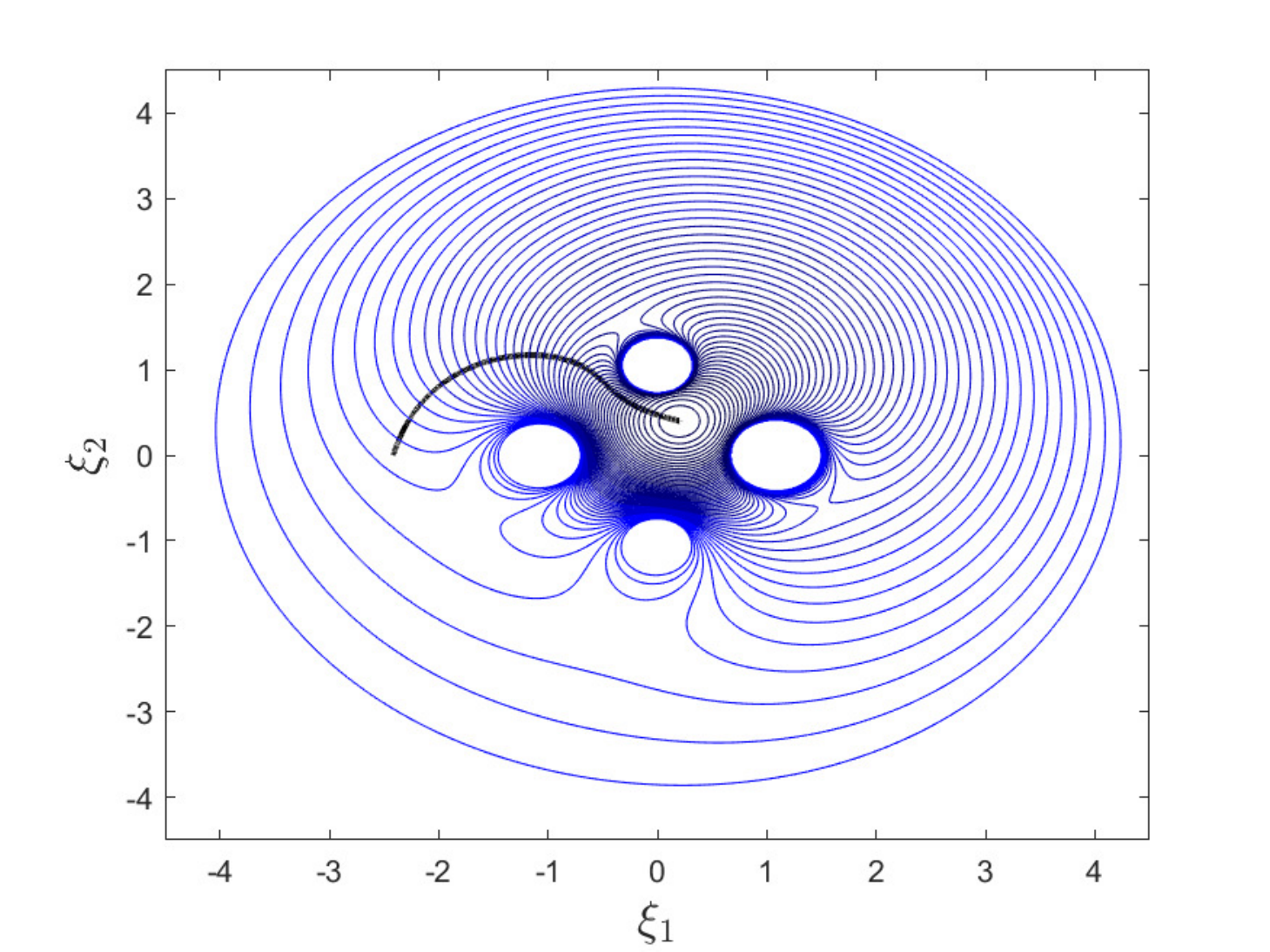}
\caption{Left: trajectory of the spherical pendulum (blue) on the sphere $\mathbb{S}^2$ in presence of constraint zones (red). Right: trajectory of the stereographic projection coordinates (black) in the corresponding Euclidean sphere world on $\mathbb{R}^2$.}
\label{fig:simulation}
\end{figure}

We consider the kinematics of the spherical pendulum
\begin{align}
    \dot x=x\times u=:\Pi(x)u
\end{align}
where $\times$ denotes the cross product and $u$ is the angular velocity of the pendulum. Using the cross product identities $\Pi(x)^\top=-\Pi(x)$ and $\Pi(x)^2=-I_3+xx^\top$, it is easy to show that $
    \Sigma(x)\Sigma(x)^\top=(1-e_3^\top x)^{-2}I_2$. It follows from \eqref{eq:nabla-psi} and \eqref{eq:Sigma-inverse} that $
    (\Sigma(x))^+=-\Pi(x)((1-x_3)I_3+e_3x^\top)J_2^\top$. 
For simulation, we pick $x(0)=(-1,0,1)/\sqrt{2}, x_d=(1,2,-2)/3$ and $\gamma=k=5$ for the control parameters. We consider $5$ constraints zones such that $a_0=e_3, a_1=e_1, a_2=-e_1, a_3=e_2$ and $a_4=-e_2$. The angles are given by $\theta_i=\pi/(7+i)$ for all $i=0,1,\cdots,4$. It is easy to check that Assumption \ref{assumption:general} holds. Simulation results are plotted in Figure \ref{fig:simulation} which show a successful constrained stabilization  on the unit $2-$sphere in the presence of different constraint zones.  The complete simulation video can be found at \url{https://youtu.be/ye8deIheiok}.
\bibliographystyle{apalike}
\bibliography{constrained-stabilization.bib}

\end{document}